\documentclass{amsart}
 \usepackage{amsmath}
 \usepackage{amssymb}

\newtheorem{theorem}{Theorem}[section]
\newtheorem{corollary}[theorem]{Corollary}
\newtheorem{lemma}[theorem]{Lemma}

\newtheorem{definition}[theorem]{Definition}
\newtheorem{conjecture}[theorem]{Conjecture}
\newtheorem{example}[theorem]{Example}

\numberwithin{equation}{section}

\title{The Multiplier Group of a Quasiperiodic Flow}

\author{Lennard F. Bakker}


\email{bakker@math.byu.edu}

\begin{document}

%\centerline{\scshape Lennard F. Bakker}
%\medskip
%{\footnotesize
 %\centerline{Department of Mathematics}
 % \centerline{Brigham Young University}
 %  \centerline{Provo, UT 84602 USA}
  % }

%\medskip

\begin{abstract} As an absolute invariant of smooth conjugacy, the multiplier group described the types of space-time symmetries that the flow has, and for a quasiperiodic flow on the $n$-torus, is the determining factor of the structure of its generalized symmetry group. It is conjectured that a quasiperiodic flow is $F$-algebraic if and only if its multiplier group is a finite index subgroup of the group of units in the ring of integers of a real algebraic number field $F$, and that a quasiperiodic flow is transcendental if and only if its multiplier group is $\{1,-1\}$. These two conjectures are partially validated for $n\geq 2$, and fully validated for $n=2$.
\end{abstract}

\maketitle

%% The name of the associate editor will be entered by a editorial staff
% \centerline{(Communicated by Aim Sciences)}
% \medskip

\section{Introduction} The multiplier group of a flow emerged during a systematic use of group (representation) theory in a recent study of space-time symmetries of flows more general than symmetries and reversing symmetries by Conner and the author in \cite{BC}, the kind of study suggested not long ago by Lamb and Roberts (p.25 \cite{RL}). As a representation of the generalized symmetry group of a flow, the multiplier group describes the types of generalized space-time symmetries that a flow possesses. The generalized symmetry group of a smooth (i.e.\ $C^\infty$) flow $\phi:{\mathbb R}\times P\to P$ on a smooth manifold $P$ without boundary is
\[ S_\phi = \{ R\in{\rm Diff}(P):{\rm there\ is\ }\alpha\in{\mathbb R}^*{\rm\ for\ which\ }R_*X_\phi=\alpha X_\phi \},\]
where ${\rm Diff}(P)$ is the group of smooth diffeomorphisms of $P$, ${\mathbb R}^*={\mathbb R}\setminus\{0\}$ is the multiplicative real group, $X_\phi(p) = (d/dt)\phi(t,p)\vert_{t=0}$ for $p\in P$
is the vector field generated by $\phi$, and $R_*X_\phi={\bf T}RX_\phi R^{-1}$ is the push-forward of $X_\phi$ by $R$ in which ${\bf T}R$ is the derivative map of $R$. A one-dimensional linear representation of $S_\phi$ is the homomorphism $\rho_\phi:S_\phi\to{\mathbb R}^*$ which takes each generalized symmetry $R\in S_\phi$ to its unique multiplier $\alpha=\rho_\phi(R)$ that appears in $R_*X_\phi=\alpha X_\phi$. The type of a generalized symmetry is determined by its multiplier: an $R\in S_\phi$ is a symmetry if $\rho_\phi(R)=1$, a reversing symmetry if $\rho_\phi(R)=-1$, and a more general type of space-time symmetry if $\rho_\phi(R)\ne\pm 1$. The multiplier group of $\phi$ is
\[ M_\phi=\rho_\phi(S_\phi).\]

The generalized symmetry group and its representation, the multiplier group, are incomplete invariants of two equivalence relations on flows on $P$. They are absolute invariants of scale equivalence in that if a flow $\phi$ on $P$ is scale equivalent to a flow $\psi$ on $P$, i.e.\ there is a $\vartheta\in{\mathbb R}^*$ such that $\vartheta X_\phi = X_\psi$, then $S_\phi=S_\psi$ and $M_\phi=M_\psi$. (The proof of this invariance is straightforward.) The multiplier group is an absolute invariant of smooth conjugacy in that if $\phi$ is smoothly conjugate to $\psi$, i.e.\ there is a $V\in{\rm Diff}(P)$ such that $V_*X_\phi = X_\psi$, then $M_\phi = M_\psi$ (Theorem 4.2 in \cite{BC}). However, the generalized symmetry group is a relative invariant of smooth conjugacy in that if $\phi$ is smoothly conjugate to $\psi$, then $S_\phi$ and $S_\psi$ are conjugate subgroups of ${\rm Diff}(P)$ (Theorem 4.1 in \cite{BC}).

Although the exact role of the multiplier group in the structure of the generalized symmetry group is not known for an arbitrary flow, it is known for a quasiperiodic flow. Recall that a flow $\phi$ on the $n$-torus, $T^n={\mathbb R}^n / {\mathbb Z}^n$, $n\geq 2$, is quasiperiodic if there is a $V\in{\rm Diff}(T^n)$ such that $V_*X_\phi$ is a constant vector field whose $n$ real components, or frequencies, are independent over ${\mathbb Q}$ (pp.79-80 in \cite{HB}). If a flow $\phi$ on $T^n$ is quasiperiodic, then there is a subgroup $H_\phi$ of $S_\phi$ isomorphic to $M_\phi$ for which $S_\phi$ is the semidirect product of the normal subgroup $\rho_\phi^{-1}(\{1\})$ by the subgroup $H_\phi$, corresponding to the conjugating homomorphism $\Gamma:H_\phi\to{\rm Aut}\big(\rho_\phi^{-1}(\{1\})\big)$ defined by $\Gamma(h)(g) = hgh^{-1}$ for $g\in\rho_\phi^{-1}(\{1\})$ and $h\in H_\phi$ (Theorem 5.5 in \cite{BA2}); in symbols this semidirect product structure is written
\[ S_\phi = \rho_\phi^{-1}(\{1\})\rtimes_\Gamma H_\phi.\]
The factor $\rho_\phi^{-1}(\{1\})$ is the group of symmetries of $\phi$, and is the group of translations on $T^n$ whenever $\phi$ is quasiperiodic, which group of translations is isomorphic to $T^n$ (Corollary 4.7 in \cite{BA2}). The multiplier group is therefore the determining factor of the structure of the generalized symmetry group of a quasiperiodic flow.

The nature of the multiplier group of a quasiperiodic flow was conjectured by the author in \cite{BA3} to be one of two possible kinds. Depending on the frequencies of a quasiperiodic flow $\phi$ on $T^n$, the multiplier group is a finitely generated group that is either isomorphic to ${\mathbb Z}_2$, or is isomorphic to ${\mathbb Z}_2\times{\mathbb Z}\times\cdot\cdot\cdot\times{\mathbb Z}$, where the number of factors of ${\mathbb Z}$ is a positive integer no bigger than $n-1$. The exact statements describing the two conjectured possible kinds of multiplier groups for quasiperiodic flows are developed in Section 2. The correction of the relationship between smooth semiconjugacy of quasiperiodic flows which was stated incorrectly in \cite{BA4} is given and followed by a partial validation of the these conjectures for quasiperiodic flows on $T^n$ for all $n\geq 2$ in Section 3. Full validation of these conjectures for quasiperiodic flows on $T^2$ is given in Section 4. Throughout the paper, examples are given to illustrate the ideas, theory, and the alluded to correction.

\section{Two Conjectures} Several necessary conditions are imposed on the multiplier group by quasiperiodicity. In view of the definition of quasiperiodicity and the absolute invariance of the multiplier group under smooth conjugacy, it is assumed without loss of generality in the sequel, that all quasiperiodic flows are generated by constant vector fields:
\[ X_\phi = \sum_{i=1}^n a_i \frac{\partial }{\partial \theta_i},\]
where $a_i\in{\mathbb R}$ and $(\theta_1,\dots,\theta_n)$ are global coordinates on $T^n$. The analysis of the multiplier group of a quasiperiodic flow reduces to problems in linear algebra since $R_*X_\phi=\alpha X_\phi$ becomes ${\bf T}RX_\phi=\alpha X_\phi$.
One necessary condition then is that
\[ M_\phi\cap{\mathbb Q}=\{1,-1\}\]
(Corollary 4.4 in \cite{BA2}). Another is that each $\alpha\in M_\phi$ is a real algebraic integer of degree at most $n$, i.e.\ a real root of a monic polynomial of degree at most $n$ in the polynomial ring ${\mathbb Z}[z]$ (Corollary 4.4 in \cite{BA2}).
A third is that for each $\alpha\in M_\phi$ there is a unique $B=(b_{ij})\in {\rm GL}(n,{\mathbb Z})$ such that
\[ \alpha = \sum_{j=1}^n b_{ij}\frac{a_j}{a_i}, {\rm\ for\ all\ }i=1,\dots,n,\]
(Corollary 4.5 in \cite{BA2}) where ${\rm GL}(n,{\mathbb Z})$ is the group of $n\times n$ matrices with integers entries and determinant $\pm 1$. This necessary condition is equivalent to ${\bf T}RX_\phi =\alpha X_\phi$ with ${\bf T}R=B$ and $\alpha=\rho_\phi(R)$, and states that each generalized symmetry of a quasiperiodic flow (generated by a constant vector field) is induced by a ${\rm GL}(n,{\mathbb Z})$ matrix plus a translation on $T^n$.

These necessary conditions are part of a connection between the group of units in the ring of integers of a real algebraic number field and the multiplier group of a quasiperiodic flow. A real algebraic number field $F$ of degree $n$ is a subfield of ${\mathbb R}$ that is a $n$-dimensional vector space over ${\mathbb Q}$. A ${\mathbb Q}$-basis for a real algebraic number field $F$ is a set of $n$ linearly independent elements of $F$. A realization of a real algebraic number field $F$ of degree $n$ over ${\mathbb Q}$ is as a simple field extension of ${\mathbb Q}$ by a real root $\delta$ of an irreducible polynomial of degree $n$ in ${\mathbb Q}[z]$. Written as $F={\mathbb Q}(\delta)$, this simple extension is the smallest subfield of ${\mathbb R}$ that contains ${\mathbb Q}\cup\{\delta\}$. A known ${\mathbb Q}$-basis for $F={\mathbb Q}(\delta)$ of degree $n$ over ${\mathbb Q}$ is $\{1,\delta,\dots,\delta^{n-1}\}$.

\begin{definition} Let $F$ be a real algebraic number field of degree $n$ over ${\mathbb Q}$. A quasiperiodic flow $\phi$ on $T^n$ is called $F$-algebraic if there is a $\vartheta\in{\mathbb R}^*$ such that $\{\vartheta a_1,\dots,\vartheta a_n\}$, the components of $\vartheta X_\phi$, form a ${\mathbb Q}$-basis for $F$.
\end{definition}

\noindent By the absolute invariance of the multiplier group under scale equivalence, the multiplier group of the flow generated by $\vartheta X_\phi$ for $\vartheta\in{\mathbb R}^*$ is the same as that of $\phi$. Associated to a real algebraic number field $F$ is its ring of integers, ${\mathfrak o}_F$, which is the set of real algebraic integers in $F$, and has the structure of a ${\mathbb Z}$-module. Associated to ${\mathfrak o}_F$ is its group of units,
\[ {\mathfrak o}_F^* = \{ \beta\in{\mathfrak o}_F\setminus\{0\}\  \vert\  \beta^{-1}\in{\mathfrak o}_F\}.\]
By Dirichlet's Unit Theorem (see p.21 in \cite{SD}), this group of units is a finitely generated group isomorphic to ${\mathbb Z}_2\times{\mathbb Z}\times\cdot\cdot\cdot\times{\mathbb Z}$ in which the number of factors of ${\mathbb Z}$ is $r_1+r_2-1$ where $r_1$ is the number of embeddings of $F$ into ${\mathbb R}$ and $2r_2$ is the number of embeddings of $F$ into ${\mathbb C}$. (Note that $r_1+r_2-1$ is no bigger than $n-1$ since $r_1+2r_2=n$, the degree of $F$ over ${\mathbb Q}$.) If a quasiperiodic flow $\phi$ on $T^n$ is $F$-algebraic, then $M_\phi\subset {\mathfrak o}_F^*$ (Theorem  3.4 in \cite{BA3}). In conjunction with the necessary condition $M_\phi\cap{\mathbb Q}=\{1,-1\}$, it follows that
\[ \{1,-1\}\subset M_\phi \subset {\mathfrak o}_F^*\]
for any quasiperiodic flow $\phi$ that is $F$-algebraic. Since ${\mathfrak o}_F^*$ is finitely generated, every subgroup of ${\mathfrak o}_F^*$ is finitely generated (Corollary 1.7, p.74 in \cite{HU}). Not every subgroup $G$ of ${\mathfrak o}_F^*$ that contains $\{1,-1\}$ has finite index, where the index of a subgroup $G$ of ${\mathfrak o}_F^*$ is the order of ${\mathfrak o}_F^*/G$ and is written $[{\mathfrak o}_F^*:G]$. The subgroup $\{1,-1\}$, which is isomorphic to ${\mathbb Z}_2$, is an infinite index subgroup of ${\mathfrak o}_F^*$. However, if $M_\phi$ is a finite index subgroup of ${\mathfrak o}_F^*$, then $M_\phi$ is isomorphic to ${\mathfrak o}_F^*$, in which case $M_\phi$ is isomorphic to ${\mathbb Z}_2\times{\mathbb Z}\times\cdot\cdot\cdot\times{\mathbb Z}$ with $r_1+r_2-1\leq n-1$ factors of ${\mathbb Z}$.

The ``upper bound'' ${\mathfrak o}_F^*$ for the multiplier group of an $F$-algebraic quasiperiodic flow $\phi$ is sharp. If the components of $X_\phi$ (or some nonzero scalar multiple of them) form a ${\mathbb Z}$-basis for ${\mathfrak o}_F$, then $M_\phi = {\mathfrak o}_F^*$ (Theorem 3.8 in \cite{BA3}). Another possibility for $M_\phi$ when $\phi$ is a quasiperiodic flow that is $F$-algebraic is illustrated next.

\begin{example}{\rm Let $\psi$ be the flow on $T^3$ with generating vector field
\[ X_\psi = \frac{\partial}{\partial\theta_1} + 3(2^{1/3})\frac{\partial}{\partial\theta_2} - 3(2^{2/3})\frac{\partial}{\partial\theta_3}.\]
Let $F={\mathbb Q}(2^{1/3})$ where $2^{1/3}$ is the real root of the irreducible $z^3-2$ in ${\mathbb Q}[z]$. One ${\mathbb Q}$-basis for $F$ is $\{1,2^{1/3},2^{2/3}\}$. Another ${\mathbb Q}$-basis for $F$ is $\{1,3(2^{1/3}),-3(2^{2/3})\}$ since this set is related to the first ${\mathbb Q}$-basis by a ${\rm GL}(3,{\mathbb Q})$ matrix. Thus $\phi$ is a quasiperiodic flow that is $F$-algebraic. The $R\in{\rm Diff}(T^3)$ induced by the ${\rm GL}(3,{\mathbb Z})$ matrix
\[ B = \begin{bmatrix} 1 & 1 & 1 \\ -18 & 1 & -3 \\ -18 & 6 & 1\end{bmatrix}\]
is in $S_\psi$ with $\rho_\psi(R)=1+3(2^{1/3})-3(2^{2/3})$. The group of units in ${\mathfrak o}_F$ is
\[ {\mathfrak o}_F^*=\{\pm(-1+2^{1/3})^k:k\in{\mathbb Z}\}\cong {\mathbb Z}_2\times{\mathbb Z} \]
(see p.201 in \cite{FT}). Notice that $\rho_\psi(R)=(-1+2^{1/3})^3$. Since by the third necessary condition on $M_\psi$, there are for each $\alpha\in M_\psi$ integers $c_{11},c_{12},c_{13}$ such that
\[ \alpha = c_{11} + 3(2^{1/3})c_{12} -3(2^{2/3})c_{13},\]
it follows that neither $-1+2^{1/3}$ nor $(-1+2^{1/3})^2=1-2(2^{1/3})+2^{2/3}$ belongs to $M_\psi$. Therefore, $M_\psi=\{\pm(-1+2^{1/3})^{3k}:k\in{\mathbb Z}\}\cong{\mathbb Z}_2\times{\mathbb Z}$ with $[{\mathfrak o}_F^*:M_\psi]=3$.
}\end{example}

This example and others (Example 4.9 in \cite{BA2} and Example 3.5 in \cite{BA3}) suggest that the $F$-algebraicity of a quasiperiodic flow $\psi$ and the finiteness of the index of $M_\psi$ as a subgroup of ${\mathfrak o}_F^*$ are closely related. Indeed, the $F$-algebraicity of $\psi$ appears to be a sufficient condition for the existence of multipliers of $\psi$ other than $\pm 1$ in the sense that $M_\psi$ is a finite index subgroup of ${\mathfrak o}_F^*$. On the other hand, the finiteness of the index of $M_\psi$ as a subgroup of ${\mathfrak o}_F^*$ appears to imply the $F$-algebraicity of $\psi$. This close relation leads to the first conjecture about the possibilities for the multiplier group of a quasiperiodic flow. 

\begin{conjecture}\label{conjFalgebraic} A quasiperiodic flow $\psi$ is $F$-algebraic if and only if $M_\psi$ is a finite index subgroup of ${\mathfrak o}_F^*$
\end{conjecture}

The collection of quasiperiodic flows on $T^n$ divides into two broad classes: those that are $F$-algebraic for some real algebraic number field $F$ of degree $n$ over ${\mathbb Q}$ (the algebraic class of quasiperiodic flows), and those that are not.

\begin{definition} A quasiperiodic flow $\psi$ on $T^n$ is called transcendental if $\psi$ is not $F$-algebraic for any real algebraic number field $F$ of degree $n$ over ${\mathbb Q}$.
\end{definition}

\noindent The second class of quasiperiodic flows on $T^n$ consists of those that are transcendental, and is called the transcendental class. The multiplier group of a quasiperiodic flow in this class always contains $\{1,-1\}$ by the first necessary condition on the multiplier group. The question is whether the multiplier group of a transcendental flow contains more than $\{1,-1\}$. A possibility for the multiplier group of a quasiperiodic flow that is transcendental is illustrated next.

\begin{example}\label{exampletrans}{\rm Let $\psi$ be the flow on $T^3$ with generating vector field
\[ X_\psi = \frac{\partial}{\partial\theta_1} + \pi \frac{\partial}{\partial\theta_2} +\pi^2\frac{\partial}{\partial\theta_3}.\]
If $\psi$ were not quasiperiodic, then there would be $u_1,u_2,u_3\in{\mathbb Q}$, not all zero, such that
\[ u_1+u_2\pi + u_3\pi^2 =0.\]
Since either $u_2\ne 0$ or $u_3\ne0$, it follows that $\pi$ would be a root of a nonconstant polynomial in ${\mathbb Q}[z]$ of degree at most $3$, and hence that $\pi$ would be algebraic. But this contradicts that $\pi$ is transcendental. So $\psi$ is quasiperiodic. Suppose that $\psi$ is $F$-algebraic for an $F={\mathbb Q}(\delta)$ of degree $3$ over ${\mathbb Q}$. Then there is a $\vartheta\in{\mathbb R}^*$ such that $\{\vartheta,\vartheta\pi,\vartheta\pi^2\}$ is a ${\mathbb Q}$-basis for $F$. Since $F$ is a field and $\vartheta\ne0$, then $\vartheta\pi/\vartheta=\pi$ is in $F$. But every element of $F$ is algebraic (Theorem 4.9, p.50 in \cite{PD}). This contradiction shows that $\psi$ is transcendental. Suppose that $M_\psi\ne\{1,-1\}$, and let $\alpha\in M_\psi\setminus\{1,-1\}$. By the third necessary condition on $M_\psi$, there are $c_{11},c_{12},c_{13}\in{\mathbb Z}$ such that
\[ \alpha = c_{11}+c_{12}\pi + c_{13}\pi^2.\]
If $c_{12}= 0$ and $c_{13}=0$, then $\alpha=\pm1$ since $M_\psi\cap{\mathbb Q}=\{1,-1\}$. Hence either $c_{12}\ne 0$ or $c_{13}\ne 0$. By the second necessary condition on $M_\psi$, there is a monic polynomial $l(z)$ in ${\mathbb Z}[z]$ of degree at most $3$ such that $l(\alpha)=0$. This means that $\pi$ is a root of a nonconstant polynomial in ${\mathbb Q}[z]$ of degree at most $6$, again a contradiction to $\pi$ being transcendental. Therefore, $M_\psi=\{1,-1\}$.
}\end{example}

This example and another (Example 2.2 in \cite{BA1}) suggest that the transcendentality of a quasiperiodic flow $\psi$ and $M_\psi=\{1,-1\}$ are closely related. Indeed, transcendentality of $\psi$ appears to prevent the existence of multipliers other than $\pm 1$. On the other hand, the lack of existence of multipliers of $\psi$ other than $\pm1$ appears to imply the transcendentality of $\psi$. It is this close relation that leads to the second conjecture about the possibilities for the multiplier group of a quasiperiodic flow.

\begin{conjecture}\label{conjtranscendental} A quasiperiodic flow $\psi$ is transcendental if and only if $M_\psi=\{1,-1\}$.
\end{conjecture}

The two conjectures together assert that the multiplier group, as an absolute invariant of smooth conjugacy, distinguishes the algebraic class of quasiperiodic flows from the transcendental class of quasiperiodic flows. They together assert a characterization of those quasiperiodic flows which have multipliers other than $\pm1$ and those which only have $\pm 1$ as multipliers, and hence they assert a complete catalogue for the determining factor in the structure of the generalized symmetry group for quasiperiodic flows.

\section{Partial Validity of Conjectures on $T^n$} It will be shown for a quasiperiodic flow $\psi$ on $T^n$, $n\geq 2$, that $F$-algebraicity of $\psi$ implies $M_\psi$ is a finite index subgroup of ${\mathfrak o}_F^*$, and that $M_\psi=\{1,-1\}$ implies the transcendentality of $\psi$. The relation of smooth semiconjugacy will play a key role in the proof. Recall that, in general, a flow $\phi$ on $P$ is smoothly semiconjugate to a flow $\psi$ on $P$ if there is a surjective $V\in C^\infty(P)$ such that ${\bf T}VX_\phi = X_\psi V$.  For a quasiperiodic flow $\phi$ generated by a constant vector field on $T^n$ to be smoothly semiconjugate to another flow $\psi$ generated by a constant vector on $T^n$ requires that the surjective $V\in C^\infty(T^n)$ satisfying ${\bf T}VX_\phi=X_\psi$ be induced by a ${\rm GL}(n,{\mathbb Q})\cap{\rm M}(n,{\mathbb Z})$ matrix plus a translation on $T^n$ (Lemma 2.4 in \cite{BA4}), where ${\rm M}(n,{\mathbb Z})$ is the set of $n\times n$ matrices with integer entries. It was originally stated by the author as Theorem 3.5 in \cite{BA4} that $M_\psi$ should be a finite index subgroup of $M_\phi$ whenever $\phi$ was an $F$-algebraic quasiperiodic flow that was smoothly semiconjugate to a flow $\psi$ on $T^n$ generated by a constant vector field. However this is incorrect as illustrated next.

\begin{example}\label{counterexample}{\rm Let $F={\mathbb Q}(\sqrt 3)$. Then ${\mathfrak o}_F^* = \{\pm(2+\sqrt 3)^k:k\in{\mathbb Z}\}$ (p.106 in \cite{HC}). Let $\phi$ and $\psi$ be $F$-algebraic quasiperiodic flows on $T^2$ with generating vector fields
\[ X_\phi = \frac{\partial}{\partial\theta_1} + 4\sqrt 3 \frac{\partial}{\partial\theta_2}\quad{\rm and}\quad X_\psi = 4\frac{\partial}{\partial\theta_1} + 60\sqrt 3\frac{\partial}{\partial\theta_2}.\]
The surjective $V\in C^\infty(T^2)$ induced by the ${\rm GL}(2,{\mathbb Q})\cap{\rm M}(2,{\mathbb Z})$ matrix
\[ \begin{bmatrix} 4 & 0 \\ 0 & 15\end{bmatrix}\]
satisfies ${\bf T}VX_\phi = X_\psi$; thus $\phi$ is smoothly semiconjugate to $\psi$. The $R\in{\rm Diff}(T^2)$ induced by the ${\rm GL}(2,{\mathbb Z})$ matrix
\[ \begin{bmatrix} 7 & 1 \\ 48 & 7\end{bmatrix}\]
satisfies $R_*X_\phi = (2+\sqrt 3)^2X_\phi$. By the third necessary condition on the multiplier group of a quasiperiodic flow, there is for each $\alpha\in M_\phi$ two integers $c_{11}$ and $c_{12}$ such that $\alpha = c_{11}+4\sqrt 3 c_{12}$. This implies that $2+\sqrt 3\not\in M_\phi$, so that $[{\mathfrak o}_F^*:M_\phi]=2$. The $Q\in{\rm Diff}(T^2)$ induced by the ${\rm GL}(2,{\mathbb Z})$ matrix
\[ \begin{bmatrix} 26 & 1 \\ 675 & 26\end{bmatrix} \]
satisfies $Q_*X_\psi = (2+\sqrt 3)^3 X_\psi$. By the third necessary condition on the multiplier group of a quasiperiodic flow, there is for each $\alpha\in M_\psi$ two integers $c_{11}$ and $c_{12}$ such that $\alpha = c_{11} + 15\sqrt 3c_{12}$. This implies that $2+\sqrt 3$ and $(2+\sqrt 3)^2 = 7+4\sqrt 3$ are not in $M_\psi$, so that $[{\mathfrak o}_F^*:M_\psi]=3$. Since $(2+\sqrt 3)^3\in M_\psi$ but $(2+\sqrt 3)^3\not\in M_\phi$, it follows that $M_\psi\not\subset M_\phi$.
}\end{example}

This counterexample shows that smooth semiconjugacy from a quasiperiodic flow $\phi$ to a flow $\psi$ generated by a constant vector field (i.e. ${\bf T}VX_\phi=X_\psi$) is not sufficient to guarantee that $M_\psi\subset M_\phi$. (Therefore, Theorem 3.1 and Corollary 3.2 in \cite{BA4} are false.) The correction of Theorem 3.5 in \cite{BA4} (along with Corollaries 3.4 and 3.6 in \cite{BA4}) is achieved by including in the hypotheses the condition that $M_\psi\subset M_\phi$.

\begin{theorem}\label{semiconjugacy} If $\phi$ is an $F$-algebraic quasiperiodic flow on $T^n$, $\psi$ is a flow on $T^n$ generated by a constant vector field, $\phi$ is smoothly semiconjugate to $\psi$, and $M_\psi\subset M_\phi$, then $\psi$ is an $F$-algebraic quasiperiodic flow and $M_\psi$ is a finite index subgroup of $M_\phi$.
\end{theorem}

\begin{proof} A key element, Theorem 3.3 in \cite{BA4}, is correct without the extra hypothesis. It states that if $V$ a smooth semiconjugacy from a quasiperiodic flow $\phi$ on $T^n$ to a flow $\psi$ generated by a constant vector field on $T^n$, then for each $R\in S_\phi$ there exists $k\in{\mathbb Z}^+$ and $Q\in S_\psi$ such that $QV=VR^k$ and $\rho_\psi(Q) = [\rho_\phi(R)]^k$. A consequence of this key element (Corollary 3.4 in \cite{BA4}) is corrected by adding $M_\psi\subset M_\psi$ to its hypotheses. Its correction states that if $\phi$ is a quasiperiodic flow on $T^n$ that is smoothly semiconjugate to a flow $\psi$ on $T^n$ generated by a constant vector field, $M_\phi$ is finitely generated, and $M_\psi\subset M_\phi$, then $M_\psi$ is a finite index subgroup of $M_\phi$. The hypothesis that the quasiperiodic flow $\phi$ be $F$-algebraic implies that $M_\phi$ is finitely generated since $M_\phi\subset{\mathfrak o}_F^*$ and ${\mathfrak o}_F^*$ is finitely generated. Thus, the addition of the extra hypothesis of $M_\psi\subset M_\phi$ to Theorem 3.5 in \cite{BA4} achieves the correction. The $F$-algebraic quasiperiodicity of $\psi$ follows directly from Theorem 2.5 in \cite{BA4}.
\end{proof}

The idea behind proving that $F$-algebraicity of a quasiperiodic flow $\psi$ implies $M_\psi$ is a finite index subgroup of ${\mathfrak o}_F^*$ is to compare a quasiperiodic flow $\phi$ with a known multiplier group via smooth semiconjugacy to a flow scale equivalent to $\psi$.

\begin{theorem}\label{finiteindex} If a quasiperiodic flow $\psi$ is $F$-algebraic, then $M_\psi$ is a finite index subgroup of ${\mathfrak o}_F^*$.
\end{theorem}

\begin{proof} By absolute invariance of the multiplier group under scale equivalence, it is assumed without loss of generality that the components of
\[ X_\psi=\sum_{i=1}^n d_i\frac{\partial}{\partial\theta_i}\]
form a ${\mathbb Q}$-basis for $F$. Since each $d_i$ is in $F$, there is a nonzero $m_i\in{\mathbb Z}$ such that $m_id_i$ is in ${\mathfrak o}_F$ (Theorem 6.5, p.77 in \cite{PD}). Set $m=m_1\cdot\cdot\cdot m_n$, and let $\Psi$ be the flow generated by
\[ X_\Psi = \sum_{i=1}^n md_i\frac{\partial}{\partial\theta_i}.\]
Then $\Psi$ is scale equivalent to $\psi$, and hence $\Psi$ is $F$-algebraic and $M_\Psi=M_\psi$. Let $\omega_1,\dots,\omega_n$ be a ${\mathbb Z}$-basis for ${\mathfrak o}_F$. Since each $md_i$ is in ${\mathfrak o}_F$, there are $b_{ij}\in{\mathbb Z}$ such that
\[ \sum_{j=1}^n b_{ij}\omega_j = md_i{\ \rm for\ }i=1,\dots,n.\]
Since the sets $\{\omega_1,\dots,\omega_n\}$ and $\{md_1,\dots,md_n\}$ are ${\mathbb Q}$-bases for $F$, the matrix $B=(b_{ij})$ is in ${\rm GL}(n,{\mathbb Q})\cap{\rm M}(n,{\mathbb Z})$. Thus the flow $\phi$ generated by
\[ X_\phi = \sum_{i=1}^n \omega_i\frac{\partial}{\partial\theta_i}\]
is smoothly semiconjugate to $\Psi$ by the surjection $V\in C^\infty(T^n)$ induced by $B$. Since $\phi$ is $F$-algebraic, $M_\phi={\mathfrak o}_F^*$ (Theorem 3.4 in \cite{BA3}), and $M_\Psi \subset{\mathfrak o}_F^*$ (Theorem 3.8 in \cite{BA3}), it follows by Theorem \ref{semiconjugacy} that $M_\psi=M_\Psi$ is a finite index subgroup of ${\mathfrak o}_F^*$.
\end{proof}

\begin{corollary}\label{structure1} If a quasiperiodic flow $\psi$ on $T^n$ is $F$-algebraic, then
\[ S_\psi \cong T^n \rtimes_\Gamma {\mathfrak o}_F^*.\]
\end{corollary}

\begin{proof}
The group $M_\psi$ is a finite index subgroup of ${\mathfrak o}_F^*$ by Theorem \ref{finiteindex}. Hence $M_\psi$ is isomorphic to ${\mathfrak o}_F^*$. In $S_\psi = \rho_\psi^{-1}(\{1\})\rtimes_\Gamma H_\psi$, the subgroup $H_\psi$ is isomorphic to $M_\psi$, and $\rho_\psi^{-1}(\{1\})$ is isomorphic to $T^n$.
\end{proof}

\begin{example}{\rm Let $F={\mathbb Q}(\delta)$ where $\delta$ is a real root of an irreducible polynomial $l(z)$ of degree $n$ in ${\mathbb Q}[z]$ for an $n\geq 2$, and let $\psi$ be the quasiperiodic flow on $T^n$ whose generating vector field is
\[ X_\psi = \sum_{i=1}^n\left(\sum_{j=1}^n c_{ij}\delta^j\right)\frac{\partial}{\partial\theta_i}\]
for $C=(c_{ij})\in{\rm GL}(n,{\mathbb Q})$. Although the components of $X_\psi$ may not form a ${\mathbb Z}$-basis for ${\mathfrak o}_F$, they form a ${\mathbb Q}$-basis for $F$, so that $\psi$ is $F$-algebraic for any $C$. Theorem \ref{finiteindex} then guarantees that $M_\psi$ is a finite index subgroup of ${\mathfrak o}_F^*$ for any choice of $C$. This gives the existence of multipliers of $\psi$ other than $\pm 1$, which multipliers in general are computationally difficult to find especially when a generating set for ${\mathfrak o}_F^*$ is not known. The flow generated by $X_\psi$ for $C=I$ (the identity matrix) and the flow generated by $X_\psi$ for $C\in{\rm GL}(n,{\mathbb Q})\setminus{\rm GL}(n,{\mathbb Z})$ are not smoothly conjugate but have multipliers groups that may be different finite index subgroups of ${\mathfrak o}_F^*$. However, by Corollary \ref{structure1}, the group structure of $S_\psi$ for any $C$ in ${\rm GL}(n,{\mathbb Q})$ is that of $T^n\rtimes_\Gamma{\mathfrak o}_F^*$, in which ${\mathfrak o}_F^*$ is isomorphic to ${\mathbb Z}_2\times{\mathbb Z}\times\cdot\cdot\cdot\times{\mathbb Z}$ for $r_1+r_2-1$ factors of ${\mathbb Z}$ where $r_1$ is the number of real roots of $l(z)$ and $r_2$ is the number of complex conjugate pairs of roots of $l(z)$.
}\end{example}

The converse of Corollary \ref{structure1} is false even if the converse of Theorem \ref{finiteindex} (and hence Conjecture \ref{conjFalgebraic}) is valid. The group structure $S_\psi\cong T^n\rtimes_\Gamma{\mathfrak o}_F^*$ for a quasiperiodic flow $\psi$ on $T^n$ does not imply that $\psi$ is $F$-algebraic because there is always a real algebraic number field $K$ of degree $n$ over ${\mathbb Q}$ which is different from $F$ such that ${\mathfrak o}_F^*\cong{\mathfrak o}_K^*$, so that $\psi$ could be $K$-algebraic instead of $F$-algebraic. However, there are no known examples in which $S_\psi\cong T^n\rtimes_\Gamma {\mathfrak o}_F^*$ and $\psi$ is transcendental. Such examples could not exist if the converse of Theorem \ref{finiteindex} were valid.

The idea behind proving that $M_\psi=\{1,-1\}$ for a quasiperiodic flow implies the transcendentality of $\psi$ is to suppose that the implication is not true and reach a contradiction with Theorem \ref{finiteindex}.

\begin{theorem}\label{transcendental} Suppose $\psi$ is a quasiperiodic flow on $T^n$. If $M_\psi=\{1,-1\}$, then $\psi$ is transcendental.
\end{theorem}

\begin{proof} Suppose that $M_\psi=\{1,-1\}$. If $\psi$ is $F$-algebraic for a real algebraic number field $F$ of degree $n$ over ${\mathbb Q}$, then $M_\psi$ has to be a finite index subgroup of ${\mathfrak o}_F^*$ by Theorem \ref{finiteindex}. But the set $\{1,-1\}$ is an infinite index subgroup of ${\mathfrak o}_F^*$. Therefore, $\psi$ is transcendental.
\end{proof}

\begin{corollary}\label{structure2} Suppose $\psi$ is a quasiperiodic flow on $T^n$. If $S_\psi\cong T^n\rtimes_\Gamma {\mathbb Z}_2$, then $\psi$ is transcendental.
\end{corollary}

\begin{proof}
In $S_\psi = \rho_\psi^{-1}(\{1\})\rtimes_\Gamma H_\psi$, the subgroup $H_\psi$ is isomorphic to $S_\psi/\rho_\psi^{-1}(\{1\})$, and hence to ${\mathbb Z}_2$. Since $M_\psi$ is isomorphic to $H_\psi$ and since $\{1,-1\}$ is the only subgroup of ${\mathbb R}^*$ isomorphic to ${\mathbb Z}_2$, it follows that $M_\psi=\{1,-1\}$. Theorem \ref{transcendental} then implies that $\psi$ is transcendental.
\end{proof}

The converse of Corollary \ref{structure2}, that $\psi$ being transcendental implies $S_\psi\cong T^n\rtimes_\Gamma {\mathbb Z}_2$, would be valid if the converse of Theorem \ref{transcendental} (and hence Conjecture \ref{conjtranscendental}) were valid. There are no known examples where $\psi$ is transcendental and $M_\psi\ne\{1,-1\}$. No such examples could exist if the converse of Theorem \ref{transcendental} were valid. There is ample evidence that the converse of Theorem \ref{transcendental} is valid, as illustrated next.

\begin{example}{\rm Let $\gamma$ be a real transcendental number, and let $\psi$ be the flow on $T^n$ for an $n\geq 2$ whose generating vector field is
\[ X_\psi = \sum_{i=1}^n\left(\sum_{j=1}^n c_{ij} \gamma^j\right)\frac{\partial}{\partial\theta_i}\]
for $C=(c_{ij})\in{\rm GL}(n,{\mathbb Q})$. The arguments used in Example \ref{exampletrans} apply to $\psi$ to show that $\psi$ is a quasiperiodic flow that is transcendental and $M_\psi = \{1,-1\}$ for any choice of $C$. Thus, $S_\psi\cong T^n\rtimes_\Gamma{\mathbb Z}_2$ for any $C$. The flow generated by $X_\psi$ for $C=I$ and the flow generated by $X_\psi$ for $C\in{\rm GL}(n,{\mathbb Q})\setminus{\rm GL}(n,{\mathbb Z})$ are not smoothly conjugate. These flows show that $M_\psi$ and $S_\psi$ are incomplete as invariants of smooth conjugacy on the collection of quasiperiodic flows.
}\end{example}

\section{Validity of Conjectures on $T^2$} It is shown for a quasiperiodic flow $\psi$ on $T^2$ that $M_\psi$ being a finite index subgroup of ${\mathfrak o}_F^*$ implies $\psi$ is $F$-algebraic, and that transcendentality of $\psi$ implies $M_\psi=\{1,-1\}$. The relationship between the existence of multipliers other than $\pm 1$ and the frequencies of $\psi$ (Theorem 2.1 in \cite{BA1}) will play a key role in the proof, and is restated here.

\begin{theorem}\label{irreducible} Suppose $X_\psi = d_1\partial/\partial\theta_1+d_2\partial/\partial\theta_2$ is the generating vector field for a quasiperiodic flow $\psi$ on $T^2$. If $M_\psi\setminus\{1,-1\}\ne\emptyset$, then $d_2/d_1$ is the root of an irreducible quadratic polynomial in ${\mathbb Q}[z]$.
\end{theorem}

The idea behind the proof that the finite index of $M_\psi$ as a subgroup of ${\mathfrak o}_F^*$ implies $\psi$ is $F$-algebraic is to use Theorem \ref{irreducible} to define a real quadratic number field $K$ for which $\psi$ is $K$-algebraic, and then to show that $F=K$. The validation of Conjecture \ref{conjFalgebraic} for quasiperiodic flows on $T^2$ then follows.

\begin{lemma}\label{Kalgebraic} If $\psi$ is a quasiperiodic flow on $T^2$ generated by
\[ X_\psi = d_1\frac{\partial}{\partial\theta_1} + d_2\frac{\partial}{\partial\theta_2}\]
for which $M_\psi\setminus\{1,-1\}\ne\emptyset$, then $\psi$ is $K$-algebraic for $K={\mathbb Q}(d_2/d_1)$.
\end{lemma}

\begin{proof} Quasiperiodicity of $\psi$ implies that $d_1\ne0$ and $d_2\ne 0$. By Theorem \ref{irreducible}, the ratio $d_2/d_1$ is a root of an irreducible quadratic polynomial in ${\mathbb Q}[z]$. Thus $K={\mathbb Q}(d_2/d_1)$ is a real quadratic number field. Since for $\vartheta=1/d_1\in{\mathbb R}^*$, the set $\{\vartheta d_1,\vartheta d_2\} = \{1,d_2/d_1\}$ is a ${\mathbb Q}$-basis for $K$, it follows that $X_\psi$ is $K$-algebraic.
\end{proof}

\begin{theorem}\label{conversefiniteindex} Suppose $\psi$ is a quasiperiodic flow on $T^2$, and $F$ is a real quadratic number field. If $M_\psi$ is a finite index subgroup of ${\mathfrak o}_F^*$, then $\psi$ is $F$-algebraic.\end{theorem}

\begin{proof} Let $F={\mathbb Q}(\delta)$ be a real quadratic number field. Let $X_\psi = d_1\partial/\partial\theta_1 + d_2\partial/\partial\theta_2$. If $M_\psi$ is a finite index subgroup of ${\mathfrak o}_F^*$, then $M_\psi\setminus\{1,-1\}\ne\emptyset$. By Lemma \ref{Kalgebraic}, the quasiperiodic flow $\psi$ is $K$-algebraic where $K={\mathbb Q}(d_2/d_1)$. Choose $\alpha\in M_\psi\setminus\{1,-1\}$. By the third necessary condition on the multiplier group of a quasiperiodic flow, there are integers $b_{11}$ and $b_{12}$ such that
\[ \alpha = b_{11} + b_{12}\frac{d_2}{d_1}.\]
Since $M_\psi\cap{\mathbb Q}=\{-1,1\}$, it follows that $b_{12}\ne 0$. On the other hand, since $M_\psi\subset{\mathfrak o}_F^*\subset F$ and since $F$ has $\{1,\delta\}$ as a ${\mathbb Q}$-basis, there are $u_1,u_2\in{\mathbb Q}$ such that
\[ \alpha = u_1 + u_2\delta.\]
Again, since $M_\psi\cap{\mathbb Q}=\{1,-1\}$, it follows that $u_2\ne0$. The quantities $d_2/d_1$ and $\delta$ are then related by
\[ u_1+u_2\delta = b_{11} + b_{12}\frac{d_2}{d_1}.\]
Since $b_{12}\ne 0$ and $u_2\ne 0$, this relationship between $d_2/d_1$ and $\delta$ implies that $d_2/d_1\in F$ and that $\delta\in K$. But $F$ is the smallest subfield of ${\mathbb R}$ that contains ${\mathbb Q}\cup\{\delta\}$, so that $F\subset K$. Also, $K$ is the smallest subfield of ${\mathbb R}$ containing ${\mathbb Q}\cup\{d_2/d_1\}$, so that $K\subset F$. Therefore $F=K$, and $\psi$ is $F$-algebraic.
\end{proof}

\begin{corollary} Suppose that $\psi$ is a quasiperiodic flow on $T^2$, and that $F$ is a real quadratic number field. Then $\psi$ is $F$-algebraic if and only if $M_\psi$ is a finite index subgroup of ${\mathfrak o}_F^*$.
\end{corollary}

\begin{proof}
Combine Theorems \ref{finiteindex} and \ref{conversefiniteindex}.
\end{proof}

The idea behind the proof that transcendentality of a quasiperiodic flow $\psi$ on $T^2$ implies $M_\psi=\{1,-1\}$ is to use Lemma \ref{Kalgebraic} to prove the contrapositive. The validation of Conjecture \ref{conjtranscendental} and the converse of Corollary \ref{structure2} on $T^2$ then follows.

\begin{theorem}\label{conversetranscendental} If a quasiperiodic flow $\psi$ on $T^2$ is transcendental, then $M_\psi=\{1,-1\}$.
\end{theorem}

\begin{proof}
Suppose that $M_\psi\ne\{1,-1\}$ for the quasiperiodic flow $\psi$ generated by $X_\psi =  d_1\partial/\partial\theta_1 + d_2\partial/\partial\theta_2$. By Lemma \ref{Kalgebraic}, the quasiperiodic flow $\psi$ is $K$-algebraic where $K={\mathbb Q}(d_2/d_1)$.
\end{proof}

\begin{corollary} If $\psi$ is a quasiperiodic flow on $T^2$, then the following are equivalent:

{\rm (a)} \ \ \hskip0.02cm $\psi$ is transcendental,

{\rm (b)} \ \  $M_\psi=\{1,-1\}$, and

{\rm (c)} \ \ \hskip0.04cm $S_\psi \cong T^2\rtimes_\Gamma {\mathbb Z}_2$.
\end{corollary}

\begin{proof}
Part (a) implies part (b) by Theorem \ref{conversetranscendental}. Part  (b) implies part (c) since $S_\psi = \rho_\psi^{-1}(\{1\})\rtimes_\Gamma H_\psi$ where $\rho_\psi^{-1}(\{1\})\cong T^2$ and $H_\psi\cong M_\psi\cong {\mathbb Z}_2$. Part (c) implies part (a) by Corollary \ref{structure2}.
\end{proof}

The validity of Conjectures \ref{conjFalgebraic} and \ref{conjtranscendental} for quasiperiodic flows on $T^2$ gives a complete catalogue of their multiplier groups which are the determining factors in the structure of the generalized symmetry groups for these flows. This validity also shows that the multiplier group as an absolute invariant of smooth conjugacy does distinguish the algebraic from the transcendental quasiperiodic flows on $T^2$.

\end{document}